\input amstex
\loadbold
\documentstyle{amsppt}
\magnification=\magstep1

\pagewidth{6.5truein}
\pageheight{9.0truein}

\long\def\ignore#1\endignore{#1}

\ignore
\input xy \xyoption{matrix} \xyoption{arrow} 
\def\edge{\ar@{-}}
\endignore

\def\ZZ{{\Bbb Z}}

\def\CC{{\Bbb C}}

\def\O{{\Cal O}}
\def\Obq{\O_{\bold q}}
\def\Oq{\O_q}
\def\Oqp{\O_{q,p}}

\def\bfq{{\bold q}}

\def\Hom{\operatorname{Hom}}
\def\kx{k^\times}
\def\chr{\operatorname{char}}

\def\spec{\operatorname{spec}}
\def\prim{\operatorname{prim}}
\def\max{\operatorname{max}}

\def\rad{\operatorname{rad}}
\def\im{\operatorname{im}}

\def\fract{\operatorname{Fract}}

\def\G{{\Cal G}}
\def\ggspec{{\Cal G}\operatorname{-spec}}
\def\ggmax{{\Cal G}\operatorname{-max}}

\def\gam{\Gamma}
\def\kgam{k\Gamma}
\def\sperp{S^{\perp}}
\def\gamp{\Gamma^+}
\def\kgamp{k\Gamma^+}
\def\swperp{S_w^{\perp}}

\def\sperpspec{\sperp\text{-}\operatorname{spec}}
\def\sperpmax{\sperp\text{-}\operatorname{max}}
\def\swperpspec{\sperp_w\text{-}\operatorname{spec}}
\def\swperpmax{\sperp_w\text{-}\operatorname{max}}
\def\gamperp{\Gamma^{\perp}}

\def\ctil{\tilde{c}}
\def\Rtil{\widetilde{R}}
\def\Atil{\widetilde{A}}
\def\Htil{\widetilde{H}}
\def\phitil{\widetilde{\phi}}
\def\Phitil{\widetilde{\Phi}}
\def\swtil{\widetilde{S}_w}
\def\swtilperp{\swtil^{\perp}}

\def\ArScTa{{\bf 1}}
\def\BroGoo{{\bf 2}}
\def\DeCKaPr{{\bf 3}}
\def\murcia{{\bf 4}}
\def\qaffine{{\bf 5}}
\def\qaffquo{{\bf 6}}
\def\Hod{{\bf 7}}
\def\HodLevone{{\bf 8}}
\def\HodLevtwo{{\bf 9}}
\def\HLT{{\bf 10}}
\def\Ing{{\bf 11}}
\def\JosA{{\bf 12}}
\def\Jbook{{\bf 13}}
\def\Kra{{\bf 14}}
\def\McPe{{\bf 15}}
\def\Van{{\bf 16}}

\topmatter

\title Quantized Primitive Ideal Spaces as Quotients of Affine Algebraic
Varieties \endtitle
\rightheadtext{Quantum spaces as quotients of affine varieties}
\author K. R. Goodearl\endauthor
\address Department of Mathematics, University of California, Santa Barbara,
CA 93106\endaddress
\email goodearl\@math.ucsb.edu\endemail


\thanks This research was partially supported by NATO Collaborative
Research Grant 960250 and National Science Foundation research grants
DMS-9622876, DMS-9970159.
\endthanks

\abstract Given an affine algebraic variety $V$ and a quantization
$\Oq(V)$ of its coordinate ring, it is conjectured that the primitive ideal
space of $\Oq(V)$ can be expressed as a topological quotient of $V$.
Evidence in favor of this conjecture is discussed, and positive solutions
for several types of varieties (obtained in joint work with E. S. Letzter)
are described. In particular, explicit topological quotient maps are given in
the case of quantum toric varieties.
\endabstract

\endtopmatter

\document

\head Introduction\endhead

A major theme in the subject of quantum groups is the philosophy that
in the passage from a classical coordinate ring to a quantized analog, the
classical geometry is replaced by structures that should be treated as
`noncommutative geometry'. Indeed, much work has been invested into the
development of theories of noncommutative differential geometry and
noncommutative algebraic geometry. We would like to pose the question
whether these theories are entirely noncommutative, or whether traces of
classical geometry are to be found in the noncommutative geometry. This
rather vague question can, of course, be focused in any number of different
directions. We discuss one particular direction here, which was developed in
joint work with E\. S\. Letzter \cite{\qaffquo}; it concerns situations in
which quantized analogs of classical varieties contain certain quotients of
these varieties.

To take an ideal-theoretic perspective on the question posed above, recall
the way in which an affine algebraic variety $V$ is captured in its
coordinate ring
$\O(V)$: the space of maximal ideals, $\max \O(V)$, is homeomorphic to $V$,
under their respective Zariski topologies. In the passage from commutative to
noncommutative algebras, the natural analog of a maximal ideal space is a
primitive ideal space -- this point of view is taken partly on practical
grounds, because noncommutative rings often have too few maximal ideals for
various purposes, but also in order to reflect the ideal theory connected
with the study of irreducible representations. (Recall that the primitive
ideals of an algebra
$A$ are precisely the annihilators of the irreducible $A$-modules.) Thus,
given a quantized version of $\O(V)$, say
$\Oq(V)$, the natural analog of $V$ is the primitive ideal space $\prim
\Oq(V)$. The aspect of the general problem which we wish to discuss here,
then, is that of finding classical geometric structure in $\prim\Oq(V)$, and
relating it to the structure of
$V$.

We report below on some success in expressing $\prim \Oq(V)$ as a quotient
of $V$ with respect to the respective Zariski topologies. Modulo a minor
technical assumption, this is done in three situations -- when $V$ is an
algebraic torus, a full affine space, or (generalizing both of those cases)
an affine toric variety. We discuss these cases in Sections 2, 3, and 4,
respectively; we sketch some of the methods and indicate relationships
among the three cases. In particular, in the last section we establish
precise formulas for topological quotient maps in the case of quantum toric
varieties, maps were only given an existence proof in \cite{\qaffquo}.

Recall that a map
$\phi : X\rightarrow Y$ between topological spaces is a {\it topological
quotient map\/} provided $\phi$ is surjective and the topology on $Y$
coincides with the quotient topology induced by $\phi$, that is, a subset
$C\subseteq Y$ is closed in $Y$ precisely when $\phi^{-1}(C)$ is closed in
$X$. In that case, $Y$ is completely determined (as a topological space) by
the topology on $X$ together with the partition of $X$ into fibres of $\phi$.
\medskip

Throughout the paper, we fix an algebraically closed base field $k$. All the
varieties we discuss will be affine algebraic varieties over $k$.

\head 1. Quantum semisimple groups\endhead

We follow the usual practice in writing `quantum semisimple groups' as an
abbreviation for `quantized coordinate rings of semisimple groups'. Thus,
`quantum $SL_n$' refers to any quantization of the coordinate ring
$\O(SL_n(k))$. We begin by displaying the primitive ideal space of the most
basic example, the standard single parameter quantization of $SL_2$.

\definition{1.1\. Example} Consider the algebra $\Oq(SL_2(k))$,
where $q\in\kx$ is a non-root of unity. This is the
$k$-algebra with generators $X_{11}$, $X_{12}$, $X_{21}$, $X_{22}$ satisfying
the following relations:
$$\gather\alignedat3 X_{11}X_{12} &= qX_{12}X_{11} &\quad\qquad X_{11}X_{21}
&= qX_{21}X_{11} &\quad\qquad X_{12}X_{21} &= X_{21}X_{12}\\
X_{21}X_{22} &= qX_{22}X_{21} &\quad\qquad X_{12}X_{22} &=
qX_{22}X_{12} \endalignedat\\
\alignedat2 X_{11}X_{22}- X_{22}X_{11} &= (q-q^{-1})X_{12}X_{21}
&\qquad\qquad X_{11}X_{22}- qX_{12}X_{21} &=1 \endalignedat\endgather$$
Since $q$ is not a root of unity, the primitive ideal space of
$\Oq(SL_2(k))$ is not large, and can be completely described as follows:

\ignore
$$\xymatrixrowsep{2pc}\xymatrixcolsep{2pc}
\xymatrix{
 \ar@{.}[r] &\langle X_{11}-\alpha,\,X_{12},\, X_{21},\,
X_{22}-\alpha^{-1} \rangle \save+<0ex,-4ex> \drop{(\alpha \in
k^\times)} \restore \ar@{.}[r] & \\
\langle X_{12} \rangle &&\langle X_{21} \rangle\\
 \ar@{.}[r] &\langle X_{12}- \beta X_{21} \rangle
\save+<0ex,-4ex> \drop{(\beta \in k^\times)} \restore
\ar@{.}[r] & 
}$$
\endignore
\smallskip

With the Zariski topology, this space seems more like a scheme than an
affine variety, since it has many non-closed points. When taken apart as
pictured, however, it can be viewed as a disjoint union of four classical
varieties: two single points, and two punctured affine lines. Taken as a
whole, $\prim \Oq(SL_2(k))$ can be related to the variety $SL_2(k)$ via the
following map:
$$\alignat2 \left( \smallmatrix \alpha&0\\0&\alpha^{-1} \endsmallmatrix
\right) &\mapsto \langle X_{11}-\alpha,\,X_{12},\, X_{21},\,
X_{22}-\alpha^{-1} \rangle &&\qquad (\alpha\in\kx)\\
\left( \smallmatrix \alpha&0\\\gamma&\alpha^{-1} \endsmallmatrix
\right) &\mapsto \langle X_{12} \rangle &&\qquad (\alpha,\gamma\in\kx)\\
\left( \smallmatrix \alpha&\beta\\0&\alpha^{-1} \endsmallmatrix
\right) &\mapsto \langle X_{21} \rangle &&\qquad (\alpha,\beta\in\kx)\\
\left( \smallmatrix \alpha&\beta\\\gamma&\delta \endsmallmatrix
\right) &\mapsto \langle X_{12}- \beta\gamma^{-1} X_{21} \rangle &&\qquad
(\alpha,\delta\in k;\ \beta,\gamma\in\kx)\\
 &&&\qquad (\alpha\delta- \beta\gamma=1)
\endalignat$$ 
We leave as an exercise for the reader to show that the map $SL_2(k)
\twoheadrightarrow \prim
\Oq(SL_2(k))$ given above is Zariski continuous; in fact, it is a topological
quotient map.
\enddefinition

\definition{1.2}
Let $G$ be a connected semisimple algebraic group over $\CC$, and let $q\in
\CC^\times$ be a non-root of unity. Quantized coordinate rings of $G$ have
been defined in both single parameter and multiparameter versions, which we
denote $\Oq(G)$ and $\Oqp(G)$ respectively. We shall not recall the
definitions here, but refer the reader to \cite{\HodLevone},
\cite{\HodLevtwo}, \cite{\HLT},
\cite{\JosA}, \cite{\Jbook}, or
\cite{\murcia}. The first classification results for primitive ideals in
this context were proved in the single parameter case:
\enddefinition

\proclaim{Theorem} {\rm [Hodges-Levasseur, Joseph]} There exists a
bijection
$$\prim \Oq(G) \longleftrightarrow \{ \text{\, symplectic leaves in\ } G
\,\},$$
 where the symplectic leaves are computed relative to the associated
Poisson structure on $G$ arising from the quantization. \endproclaim

\demo{Proof} This was proved for the cases $G=SL_3(\CC)$ and $G= SL_n(\CC)$
by Hodges and Le\-vas\-seur in \cite{\HodLevone, Theorem 4.4.1} and
\cite{\HodLevtwo, Theorem 4.2}, and then for arbitrary $G$ by Joseph
\cite{\JosA, Theorem 9.2},
\cite{\Jbook, Theorems 10.3.7, 10.3.8}.
\qed\enddemo

Let us rewrite this result to bring in the points of $G$ more directly:
There exists a surjection $G \twoheadrightarrow \prim \Oq(G)$ whose fibres
are precisely the symplectic leaves in $G$. In the multiparameter case,
the situation is slightly more complicated, as follows.

\proclaim{1.3\. Theorem} {\rm [Hodges-Levasseur-Toro]} There exists a
surjection 
$G \twoheadrightarrow \prim \Oqp(G)$, but the fibres are symplectic leaves
only for certain choices of $p$. \endproclaim

\demo{Proof} The first statement follows from the results in \cite{\HLT,
Section 4}; for the second, see \cite{\HLT, Theorems 1.8, 4.18}.
\qed\enddemo

\definition{1.4} We can fill in a bit more detail about these surjections by
bringing in a group of symmetries. Let $H$ be a maximal torus in $G$, acting
on
$\Oqp(G)$ by `winding automorphisms' (cf\. \cite{\Jbook, (1.3.4)}). For
example, if $G= SL_2(\CC)$ and $H= \bigl\{ \left( \smallmatrix \alpha&0\\
0&\alpha^{-1} \endsmallmatrix \right) \bigm| \alpha\in \CC^\times \bigr\}$,
then $H$ acts on $\Oq(G)$ so that $\left( \smallmatrix \alpha&0\\
0&\alpha^{-1} \endsmallmatrix \right) {.} \left( \smallmatrix
X_{11}&X_{12}\\ X_{21}&X_{22} \endsmallmatrix \right)= \left( \smallmatrix
\alpha X_{11}&\alpha X_{12}\\ \alpha^{-1} X_{21}&\alpha^{-1} X_{22}
\endsmallmatrix
\right)$. In this case, the induced action of $H$ on $\prim \O_q(G)$ has
four orbits, which are the sets displayed in Example 1.1.

It
follows from the analyses in \cite{\HodLevone, \HodLevtwo, \JosA, \Jbook,
\HLT} that

$\bullet$ The number of $H$-orbits in $\prim \Oqp(G)$ is finite.

$\bullet$ Each of these $H$-orbits is Zariski homeomorphic to $H$
modulo the relevant stabilizer.

$\bullet$ The preimage of each $H$-orbit, under the surjection in Theorem
1.3, is a locally closed subset of $G$.

These properties hint strongly that there is some topological connection
between $G$ and $\prim \Oqp(G)$. 
\enddefinition

\definition{1.5\. Conjecture} There exists an $H$-equivariant
surjection $G \twoheadrightarrow \prim \Oqp(G)$ which is a topological
quotient map (with respect to the Zariski topologies). \enddefinition

If this conjecture were verified, then it together with knowledge of the
fibres of the map would completely determine $\prim \Oqp(G)$ as a
topological space.

In the case of $\Oq(SL_2(\CC))$, the conjecture can be readily verified `by
hand', as already suggested above. There is a further piece of evidence
 in the $SL_3$ case: Brown and the author have constructed an
$H$-equivariant topological quotient map $B^+ \twoheadrightarrow \prim
\Oq(B^+)$ where $B^+$ is the upper triangular Borel subgroup of $SL_3(\CC)$
\cite{Work in progress}.

The major difficulty in the way of the conjecture is that no good global
description of the Zariski topology on $\prim \Oqp(G)$ is known, only the
Zariski topologies on the separate $H$-orbits. 

In order to try to obtain a better feel for the problem, let us widen the
focus:

\definition{1.6\. Conjecture} If $\Obq(V)$ is a quantized coordinate
ring of an affine algebraic variety (or group) $V$, there exists a
topological quotient map $V \twoheadrightarrow \prim\Obq(V)$, equivariant
with respect to an appropriate group of automorphisms.
\enddefinition

Here $\Obq(V)$ is meant to stand for any quantization of $\O(V)$, a phrase
which is only slightly less vague than the symbolism, since no definition of
a ``quantization'' of a coordinate ring exists. For descriptions of many
commonly studied classes of examples, see \cite{\murcia}.

In the remainder of the paper, we discuss some cases where the
extended conjecture has been verified.

\head 2. Quantum tori \endhead

The simplest case in which to examine our conjecture is that of a
quantization of the coordinate ring of an algebraic torus. The parameters
required can be arranged as a matrix $\bfq= (q_{ij}) \in M_n(\kx)$ which is
{\it multiplicatively antisymmetric\/}, that is, $q_{ii}=1$ and $q_{ji}=
q_{ij}^{-1}$ for all $i,j$. (These conditions are assumed in order to avoid
degeneracies in the resulting algebra.) The {\it quantum torus\/} over $k$
with respect to these parameters is the algebra
$$\Obq((\kx)^n) := k\langle x_1^{\pm1}, \dots, x_n^{\pm1} \mid x_ix_j=
q_{ij}x_jx_i \text{\ for all\ } i,j \rangle.$$
The structure of this algebra has been analyzed in various cases by many
people, among them McConnell and Pettit \cite{\McPe}, De Concini, Kac, and
Procesi \cite{\DeCKaPr}, Hodges \cite{\Hod}, Vancliff \cite{\Van}, Brown and
the author \cite{\BroGoo}, Letzter and the author \cite{\qaffine}, and
Ingalls \cite{\Ing}. In particular, it is known that  $\prim \Obq((\kx)^n)$
and
$\spec \Obq((\kx)^n)$ are homeomorphic to the maximal and prime spectra,
respectively, of the center $Z(\Obq((\kx)^n))$, via contraction and
extension, and that the center is a Laurent polynomial ring over $k$ (e.g.,
\cite{\qaffine, Lemma 1.2, Corollary 1.5}).

\proclaim{2.1\. Theorem} {\rm \cite{\qaffquo, Theorem 3.11}} For any
multiplicatively antisymmetric matrix $\bfq\in M_n(\kx)$, there exist
topological quotient maps
$$\align (\kx)^n &\twoheadrightarrow
\prim \Obq((\kx)^n)\\
 \spec k[y_1^{\pm1}, \dots, y_n^{\pm1}] &\twoheadrightarrow
\spec \Obq((\kx)^n),\endalign$$
equivariant with respect to the natural actions of the
torus $(\kx)^n$.
\qed\endproclaim

This case, in many ways, works out too smoothly to really put Conjecture 1.6
to the test. In that regard, the case of quantum affine spaces is much more
interesting. Before moving on to that case, however, we describe the form of
the maps in Theorem 2.1, since this will preview essential aspects of the
quantum affine space case. We proceed by rewriting quantum tori as twisted
group algebras, as follows.

\definition{2.2} Fix a quantum torus $A= \Obq((\kx)^n)$. Set $\gam= \ZZ^n$,
and define $\sigma: \gam\times\gam \rightarrow \kx$ by the rule
$$\sigma(\alpha,\beta)= \prod_{i,j=1}^n q_{ij}^{\alpha_i\beta_j}.$$
Then $\sigma$ is an alternating bicharacter on $\gam$, and it determines the
commutation rules in $A$, since $x^\alpha x^\beta= \sigma(\alpha,\beta)
x^\beta x^\alpha$ for all $\alpha,\beta \in\gam$, where we have used the
standard multi-index notation $x^\alpha= x_1^{\alpha_1} x_2^{\alpha_2}
\cdots x_n^{\alpha_n}$. Recall that the {\it radical\/} of $\sigma$ is the
subgroup
$$\rad(\sigma)= \{ \alpha\in \gam \mid \sigma(\alpha,-) \equiv 1\}$$
of $\gam$. It plays the following role: $Z(A)= k[ x^\alpha \mid \alpha\in
\rad(\sigma) ]$ (cf\. \cite{\qaffine, Lemma 1.2} or \cite{\Ing, Proposition
5.2}).

Now choose a 2-cocycle $c: \gam\times\gam \rightarrow \kx$ such that
$c(\alpha,\beta) c(\beta,\alpha)^{-1}= \sigma(\alpha,\beta)$ for
$\alpha,\beta \in \gam$. (This can be done in many ways -- see \cite{\ArScTa,
Proposition 1, p\. 888} or \cite{\qaffquo, (3.2)}, for instance.) We may
identify
$A$ with the twisted group algebra
$k^c\gam$, that is, the $k$-algebra with a basis $\{x_\alpha \mid \alpha\in
\gam\}$ such that
$$x_\alpha x_\beta= c(\alpha,\beta) x_{\alpha+\beta}$$
for $\alpha,\beta \in\gam$. 

Let us write the group algebra $\kgam$ in terms of a basis $\{y_\alpha \mid
\alpha\in
\gam\}$, where
$y_\alpha y_\beta= y_{\alpha+\beta}$
for $\alpha,\beta \in\gam$. There is a $k$-linear (vector space) isomorphism
$$\Phi_c : A\rightarrow \kgam$$
such that $\Phi_c(x_\alpha)= y_\alpha$ for $\alpha \in\gam$. While this map
is not multiplicative, it satisfies $\Phi_c(x_\alpha x_\beta)=
c(\alpha,\beta) \Phi_c(x_\alpha) \Phi_c(x_\beta)$ for $\alpha,\beta \in\gam$
(cf\. \cite{\qaffquo, (3.5)}).

The final ingredient needed to describe the topological quotient maps in the
current situation is the torus $H= (\kx)^n$ and its natural actions on $A$
and $\kgam$ by $k$-algebra automorphisms. The induced
actions of
$H$ on the prime and primitive spectra of $A$ and $\kgam$ will also be
required. We identify
$H$ with
$\Hom(\gam,\kx)$, which allows us to write the above actions as
$$h{.}x_\alpha= h(\alpha)x_\alpha \qquad\text{and}\qquad h{.}y_\alpha=
h(\alpha)y_\alpha$$
for $h\in H$ and $\alpha \in\gam$. Observe that the map $\Phi_c$ is
$H$-equivariant. The identification
$H=
\Hom(\gam,\kx)$ provides a pairing $H\times\gam \rightarrow \kx$, and we use
$^\perp$ to denote orthogonals with respect to this pairing. In particular,
$$S^\perp= \{h\in H\mid \ker h\supseteq S\}$$
for $S\subseteq \gam$. 
\enddefinition

\definition{2.3} It is convenient to have a compact notation for the
intersection of the orbit of an ideal under a group of automorphisms. If $P$
is an ideal of a ring $B$, and $T$ is a group acting on $B$ by ring
automorphisms, we write
$$(P:T)= \bigcap_{t\in T} t(P).$$
This ideal can also be described as the largest $T$-invariant ideal of $B$
contained in $P$. The same notation can also be applied to any subset of $B$.
\enddefinition

\proclaim{2.4\. Theorem} {\rm \cite{\qaffquo, Theorem 3.11}} Let $A=
\Obq((\kx)^n)$, and fix $\gam,\sigma,c,\Phi_c,H$ as above. Assume that
$c\equiv 1$ on
$\rad(\sigma)\times
\gam$. (Such a choice for $c$ is always possible; cf\. \cite{\qaffquo, Lemma
3.12}.) Then the rule $P
\mapsto
\Phi_c^{-1}(P:\rad(\sigma)^\perp)$, for prime ideals $P$ of $\kgam$, defines
$H$-equivariant topological quotient maps
 $$\spec\kgam
\twoheadrightarrow \spec A \qquad\text{and}\qquad \max\kgam
\twoheadrightarrow \prim A.$$
The fibres of the second map are exactly the $\rad(\sigma)^\perp$-orbits in
$\max\kgam$. \qed\endproclaim

\definition{2.5} As the formula in Theorem 2.4 indicates, the given maps are
compositions of two maps, which themselves have nice properties. To express
this in more detail, write $S= \rad(\sigma)$, and recall that an {\it
$\sperp$-prime ideal\/} of $\kgam$ is any proper $\sperp$-invariant ideal
$Q$ such that whenever $Q$ contains a product of $\sperp$-invariant ideals,
it must contain one of the factors. Denote by $\sperpspec \kgam$ the set of
$\sperp$-prime ideals of $\kgam$, and observe that this set supports a
Zariski topology (defined in the obvious way). The first map in Theorem 2.4
can be factored in the form
$$\spec\kgam \twoheadrightarrow \sperpspec \kgam @>{\approx}>> \spec A,$$
where the map $P\mapsto (P:\sperp)$ from $\spec\kgam$ to $\sperpspec \kgam$
is a topological quotient map, and the map $Q\mapsto \Phi_c^{-1}(Q)$ from
$\sperpspec \kgam$ to $\spec A$ is a homeomorphism. That $P\mapsto
(P:\sperp)$ is a topological quotient map is an easy exercise (cf\.
\cite{\qaffquo, Proposition 1.7} for a generalization); the bulk of the work
in proving Theorem 2.4 goes into establishing the homeomorphism $\sperpspec
\kgam \approx \spec A$ (cf\. \cite{\qaffquo, Proposition 3.10}). The key
point is that the hypothesis on $c$ ensures that $\Phi_c$ is
`central-semilinear', that is, $\Phi_c(za)= \Phi_c(z)\Phi_c(a)$ for $z\in
Z(A)$ and $a\in A$ \cite{\qaffquo, Lemma 3.5}.

The factorization above respects maximal and primitive ideals in the
following way. Let $\sperpmax \kgam$ denote the set of maximal proper
$\sperp$-invariant ideals of $\kgam$; this is a subset of $\sperpspec
\kgam$, which we equip with the relative topology. Then we have the
factorization
$$\max\kgam \twoheadrightarrow \sperpmax \kgam @>{\approx}>> \prim A,$$
where $P\mapsto (P:\sperp)$ gives a topological quotient map from
$\max\kgam$ to $\sperpmax \kgam$, and $Q\mapsto \Phi_c^{-1}(Q)$ gives a
homeomorphism from $\sperpmax \kgam$ to $\prim A$ (cf\. \cite{\qaffquo,
Propositions 3.9, 3.10}).
\enddefinition

\head 3. Quantum affine spaces\endhead

Let $\bfq\in M_n(\kx)$ again be a multiplicatively antisymmetric matrix. The
corresponding {\it quantum affine space\/} over $k$ is the algebra
$$\Obq(\kx) := k\langle x_1, \dots, x_n \mid x_ix_j=
q_{ij}x_jx_i \text{\ for all\ } i,j \rangle.$$

\definition{3.1} Denote the algebra $\Obq(\kx)$ by $A$ for purposes of
discussion. The torus $H= (\kx)^n$ again acts naturally as $k$-algebra
automorphisms on
$A$. The primitive spectrum of $A$ consists of $2^n$ $H$-orbits, which we
denote $\prim_w A$, indexed by the subsets $w\subseteq \{1,\dots,n\}$
\cite{\qaffine, Theorem 2.3}. Here
$$\prim_w A= \{ P\in\prim A\mid x_i\in P \Longleftrightarrow i\in w\}.$$ 
Similarly, $\spec A$ is a disjoint union of $2^n$ subsets
$$\spec_w A= \{ P\in\spec A\mid x_i\in P \Longleftrightarrow i\in w\}.$$ 
Each $\spec_w A$ is homeomorphic, via localization, to the prime spectrum of
the quantum torus
$$A_w= \bigl( A/ \langle x_i\mid i\in w\rangle \bigr) [x_j^{-1}\mid j\notin
w],$$
and likewise $\prim_w A$ is homeomorphic to $\prim A_w$ (cf\.
\cite{\qaffine, Theorem 2.3}).

The partition $\prim A= \bigsqcup_w \prim_w A$ is a quantum analog of the
standard stratification of affine $n$-space by its $H$-orbits. To emphasize
the parallel, we label the $H$-orbits in $k^n$ as
$$(k^n)_w= \{ (\alpha_1,\dots,\alpha_n)\in k^n \mid \alpha_i=0
\Longleftrightarrow i\in w\},$$
for $w\subseteq \{1,\dots,n\}$. It follows from Theorem 2.1 that there are
$H$-equivariant topological quotient maps
$$(k^n)_w \twoheadrightarrow \prim A_w \approx \prim_w A$$
for all $w$. The problem is then to patch these maps together coherently.
Modulo a minor technical assumption, that can be done:
\enddefinition

\proclaim{3.2\. Theorem} {\rm \cite{\qaffquo, Theorem 4.11}} Let $\bfq=
(q_{ij})\in M_n(\kx)$ be a multiplicatively antisymmetric matrix. Assume that
either
$-1\notin \langle q_{ij}\rangle$ or $\chr k=2$. Then there exist topological
quotient maps
$$\align k^n &\twoheadrightarrow
\prim \Obq(k^n)\\
 \spec k[y_1, \dots, y_n] &\twoheadrightarrow
\spec \Obq(k^n),\endalign$$
equivariant with respect to the natural actions of the
torus $(\kx)^n$.
\qed\endproclaim

Here $\langle q_{ij}\rangle$ denotes the subgroup of $\kx$ generated by the
$q_{ij}$.

 To describe the maps in the above theorem, and identify the fibres of the
first, we use a setup analogous to that in the quantum torus case. In
particular, we write
$\Obq(k^n)$ as a twisted semigroup algebra with respect to a suitable
cocycle.

\definition{3.3} Fix a quantum affine space $A= \Obq(k^n)$. Set $\gam=
\ZZ^n$ and $\gamp= (\ZZ^+)^n$, and define $\sigma: \gam\times\gam \rightarrow
\kx$ as in (2.2). As in the previous case, $\sigma$ determines the
commutation rules in
$A$.

Next, choose a 2-cocycle $c: \gam\times\gam \rightarrow \kx$ such that
$c(\alpha,\beta) c(\beta,\alpha)^{-1}= \sigma(\alpha,\beta)$ for
$\alpha,\beta \in \gam$. We identify $A$ with the twisted semigroup algebra
$k^c\gamp$ with a basis $\{x_\alpha \mid \alpha\in
\gamp\}$, and we write the corresponding semigroup algebra $R= \kgamp$ in
terms of a basis $\{y_\alpha \mid
\alpha\in
\gamp\}$. There are partitions $\spec R= \bigsqcup_w \spec_w R$ and $\max R=
\bigsqcup_w \max_w R$ analogous to those for $\spec A$ and $\prim A$.

Again as before, identify the torus $H= (\kx)^n$ with
$\Hom(\gam,\kx)$, so that the natural actions of $H$ on $A$ and $\kgamp$ are
given by
$$h{.}x_\alpha= h(\alpha)x_\alpha \qquad\text{and}\qquad h{.}y_\alpha=
h(\alpha)y_\alpha.$$
\enddefinition

\definition{3.4} It is convenient to identify the localizations $A_w$
with subalgebras of the twisted group algebra $k^c\gam$, and to identify the
corresponding localizations $R_w$ of $R$ with subalgebras of the group
algebra $\kgam$. This is done as follows:
$$A_w = \sum_{\alpha\in \gam_w} kx_\alpha \qquad\text{and}\qquad
R_w = \sum_{\alpha\in \gam_w} ky_\alpha= \kgam_w,$$
where $\gam_w= \{\alpha\in \gam \mid
\alpha_i=0 \text{\ for\ } i\in w\}$. (Here we have extended the $x_\alpha$
and the $y_\alpha$ to $\gam$-indexed bases for $k^c\gam$ and $\kgam$.) As in
(2.2), there is an
$H$-equivariant
$k$-linear  isomorphism
$\Phi= \Phi_c : k^c\gam\rightarrow \kgam$
such that $\Phi(x_\alpha)= y_\alpha$ for $\alpha \in\gam$. This map
restricts to $H$-equivariant $k$-linear isomorphisms from $A$ onto $R$
and from $A_w$ onto $R_w$ for each $w$.

Next, let $\sigma_w$ and $c_w$ denote the restrictions of $\sigma$ and $c$
to $\gam_w$. Then the identification above can be restated as $A_w=
k^{c_w}\gam_w$. We shall need the radical of $\sigma_w$, and we emphasize
that this is a subgroup of $\gam_w$. There is a natural pairing between
$\gam_w$ and a quotient group of $H$, but it is more useful to take
orthogonals with respect to the pairing of $H$ and $\gam$; thus, we write
$$\rad(\sigma_w)^\perp= \{h\in H\mid \ker h\supseteq \rad(\sigma_w)\}.$$
\enddefinition

\proclaim{3.5\. Theorem} {\rm \cite{\qaffquo, Theorem 4.11}} Let $A=
\Obq(k^n)$, and assume that either
$-1\notin \langle q_{ij}\rangle$ or $\chr k=2$. Fix
$\gam,\sigma,c,\Phi,H,\sigma_w$ as above. Assume that
$c$ is an alternating bicharacter on $\gam$ with
$c^2=\sigma$ and such that $c(\alpha,\beta)=1$ whenever
$\sigma(\alpha,\beta)=1$. (Such a choice for
$c$ is possible by
\cite{\qaffquo, Lemma 4.2}.) Then there are $H$-equivariant topological
quotient maps
 $$\spec\kgamp
\twoheadrightarrow \spec A \qquad\text{and}\qquad \max\kgamp
\twoheadrightarrow \prim A$$
such that $P
\mapsto
\Phi^{-1}(P:\rad(\sigma_w)^\perp)$ for $P\in \spec_w \kgamp$.
The fibres of the second map over points in $\prim_w A$ are exactly the
$\rad(\sigma_w)^\perp$-orbits within
$\max_w\kgamp$. \qed\endproclaim

While the precise formula given for the maps in Theorem 3.5 is particularly
useful for computations, it may round out the picture to give a global
formula (independent of $w$), as follows.

\proclaim{3.6\. Lemma} Let $\phi : \spec \kgamp \twoheadrightarrow \spec A$
be the topological quotient map given in Theorem 3.5. For all $P\in \spec
\kgamp$, the prime ideal $\phi(P)\in \spec A$ equals the largest ideal of
$A$ contained in the set $\Phi^{-1}(P)$. \endproclaim

\demo{Proof} We shall need the generators $y_i= y_{\epsilon_i}\in \kgamp$
and $x_i= x_{\epsilon_i}\in A$, where $\epsilon_1,\dots,\epsilon_n$ is the
standard basis for $\gam$.

Let $P\in \spec \kgamp$; then  $P\in \spec_w \kgamp$ for some $w$. First
suppose that $P\in \max_w \kgamp$. Then $\phi(P)\in \prim_w A$, and so
$\phi(P)$ is a maximal element of $\spec_w A$ \cite{\qaffine, Theorem 2.3}.
Now
$\phi(P)$ contains the ideal $J_w := \langle x_i\mid i\in w\rangle$. On the
other hand, the multiplicative set generated by $\{y_j\mid j\notin w\}$ in
$\kgamp$ is disjoint from $P$, and $\Phi^{-1}$ sends elements of this set to
scalar multiples of elements in the multiplicative set $X_w$ generated by
$\{x_j \mid j\notin w\}$ in $A$. Hence, $\Phi^{-1}(P)$ is disjoint from
$X_w$.

Let $I$ be the largest ideal of $A$ which is contained in $\Phi^{-1}(P)$.
Then $J_w\subseteq \phi(P) \subseteq I$ and $I$ is disjoint from $X_w$.
Hence, $I/J_w$ induces a proper ideal, call it $IA_w$, in the localization
$A_w= (A/J_w)[X_w^{-1}]$. Let $M$ be a maximal ideal of $A_w$ containing
$IA_w$, and let $Q$ be the inverse image of $M$ under the localization map
$A\rightarrow A/J_w \rightarrow A_w$. Then $Q\in \spec_w A$ and $\phi(P)
\subseteq I\subseteq Q$. Since $\phi(P)$ is maximal in $\spec_w A$, we
obtain $\phi(P)=Q$ and thus $\phi(P)=I$. Therefore the lemma holds for $P\in
\max_w \kgamp$.

Now consider an arbitrary $P\in \spec_w \kgamp$. Since $\kgamp$ is a
commutative affine $k$-algebra, $P$ is an intersection of maximal ideals. It
follows easily, as in
\cite{\BroGoo, Proposition 1.3(a)}, that
$P= \bigcap_r P_r$ for some maximal ideals $P_r\in \max_w \kgamp$. (Namely,
$P= \bigcap_{v\subseteq \{1,\dots,n\}} Q_v$ where each $Q_v$ is an
intersection of maximal ideals from $\max_v \kgamp$. Observe that $Q_v=
\kgamp$ when $v\not\supseteq w$, and that $Q_v \supsetneq P$ when
$v\supsetneq w$. Consequently, $P=Q_w$.) Since the set functions
$(-:
\rad(\sigma_w)^\perp)$ and
$\Phi^{-1}$ preserve intersections, we see that $\phi(P)= \bigcap_r
\phi(P_r)$ and $\Phi^{-1}(P)= \bigcap_r
\Phi^{-1}(P_r)$. By the previous paragraph, each $\phi(P_r)$ is the largest
ideal of $A$ contained in $\Phi^{-1}(P_r)$. Therefore $\phi(P)$ is the
largest ideal of $A$ contained in $\Phi^{-1}(P)$. \qed\enddemo

\definition{3.7} The maps in Theorem 3.5 have analogous factorizations to
those in (2.5), which we write
$$\align \spec \kgamp &\twoheadrightarrow \ggspec \kgamp @>{\approx}>> \spec
A\\
\max \kgamp &\twoheadrightarrow \ggmax \kgamp @>{\approx}>> \prim A.
\endalign$$
Here $\G$ stands for the indexed family of groups $\{ \swperp \mid
w\subseteq \{1,\dots,n\} \}$ where $S_w= \rad(\sigma_w)$, while
$$\align \ggspec \kgamp &= \bigsqcup_{w\subseteq \{1,\dots,n\}} \{
(P:\swperp) \mid P\in \spec_w \kgamp \}\\
\ggmax \kgamp &= \bigsqcup_{w\subseteq \{1,\dots,n\}} \{ (P:\swperp) \mid
P\in \max_w \kgamp \}. \endalign$$
Closed sets for Zariski topologies on $\ggspec \kgamp$ and $\ggmax \kgamp$
are defined in the usual way; a compatibility condition on the groups in
$\G$ ensures that the result actually is a topology \cite{\qaffquo,
(2.4), (4.9)}. A fairly general piece of commutative algebra, developed in
\cite{\qaffquo, Section 2}, provides the topological quotient maps from
$\spec \kgamp$ and
$\max \kgamp$ onto $\ggspec \kgamp$ and $\ggmax \kgamp$. The homeomorphisms
of $\ggspec \kgamp$ and $\ggmax \kgamp$ onto $\spec A$ and $\prim A$ are
constructed by patching together homeomorphisms of $\swperpspec \kgam_w$
and $\swperpmax \kgam_w$ onto $\spec A_w$ and $\prim A_w$ obtained as in
(2.5).
\enddefinition

\definition{3.8\. Example} To illustrate the form of these topological
quotients, we give the basic example in dimension 3, where various
differences between the classical and quantum cases appear. Recall that the
single parameter quantum affine spaces, relative to a scalar $q\in\kx$, are
the algebras
$$\Oq(k^n)= k\langle x_1,\dots,x_n \mid x_ix_j= qx_jx_i \text{\ for\ }
i<j\}.$$

{\bf (a)} Assume first that $q$ is not a root of unity, and let $p$ be a
square root of $q$ in $k$. In this case, the topological quotient map
$$k^3 \approx \max k(\ZZ^+)^3 \twoheadrightarrow \prim \Oq(k^3)$$
described in Theorem 3.5 can be computed as follows \cite{\qaffquo, (5.2)},
where the scalars $\lambda_i$ are all assumed to be nonzero:
$$\alignat2 (0,0,0) &\mapsto \langle x_1,\, x_2,\, x_3 \rangle &\qquad
(\lambda_1,\lambda_2,0) &\mapsto \langle x_3 \rangle\\
(\lambda_1,0,0) &\mapsto \langle x_1-\lambda_1,\, x_2,\, x_3 \rangle &\qquad
(\lambda_1,0,\lambda_3) &\mapsto \langle x_2 \rangle\\
(0,\lambda_2,0) &\mapsto \langle x_1,\, x_2-\lambda_2,\, x_3 \rangle &\qquad
(0,\lambda_2,\lambda_3) &\mapsto \langle x_1 \rangle \\
(0,0,\lambda_3) &\mapsto \langle x_1,\, x_2,\, x_3-\lambda_3 \rangle &\qquad
(\lambda_1,\lambda_2,\lambda_3) &\mapsto \langle \lambda_2x_1x_3-
p\lambda_1\lambda_3x_2 \rangle \endalignat$$
The right column demonstrates both the compression effect of the quotient
process and the fact that $\Oq(k^3)$ has `many fewer' primitive ideals than
$\O(k^3)$. Note the factor
$p$ in the last term. Without this factor, we obtain only a surjective map,
which is not Zariski continuous.

{\bf (b)} For contrast, consider the case where $q$ is a primitive $t$-th
root of unity for some odd $t>1$, and take $p= q^{(t+1)/2}$. Then
the topological quotient map $k^3 \twoheadrightarrow \prim \Oq(k^3)$ has the
form below \cite{\qaffquo, (5.3)}:
$$\gather \alignedat2 (0,0,0) &\mapsto \langle x_1,\, x_2,\, x_3 \rangle
&\qquad (\lambda_1,\lambda_2,0) &\mapsto \langle x_1^t-\lambda_1^t,\,
x_2^t-\lambda_2^t,\, x_3 \rangle\\
(\lambda_1,0,0) &\mapsto \langle x_1-\lambda_1,\, x_2,\, x_3 \rangle &\qquad
(\lambda_1,0,\lambda_3) &\mapsto \langle x_1^t-\lambda_1^t,\, x_2,\,
x_3^t-\lambda_3^t \rangle\\
(0,\lambda_2,0) &\mapsto \langle x_1,\, x_2-\lambda_2,\, x_3 \rangle &\qquad
(0,\lambda_2,\lambda_3) &\mapsto \langle x_1,\, x_2^t-\lambda_2^t,\,
x_3^t-\lambda_3^t \rangle \\
(0,0,\lambda_3) &\mapsto \langle x_1,\, x_2,\, x_3-\lambda_3 \rangle
&\qquad \endalignedat\\
(\lambda_1,\lambda_2,\lambda_3) \mapsto \langle x_1^t-\lambda_1^t,\,
x_2^t-\lambda_2^t,\,x_3^t-\lambda_3^t,\,  \lambda_2x_1x_3-
p\lambda_1\lambda_3x_2 \rangle \endgather$$ 
In this case, $\prim \Oq(k^3)$ is more `classical', in that all its points
are closed.
\enddefinition

\head 4. Quantum toric varieties and cocycle twists\endhead

In this section, we discuss quantizations of toric varieties \cite{\Ing} and
extend Theorem 3.5 to that setting.
There are several ways in which (affine) toric varieties may be defined.
For our purposes, the simplest definition is a variety on which a torus
acts (morphically) with a dense orbit. Some equivalent conditions are given
in the following theorem. 
 
\proclaim{4.1\. Theorem} {\rm \cite{\Kra, Satz  5, p\. 105}} Let $H$ be an
algebraic torus, acting morphically on an irreducible affine variety $V$.
The following conditions are equivalent:

{\rm (a)} There are only finitely many $H$-orbits in $V$.

{\rm (b)} Some $H$-orbit is dense in $V$.

{\rm (c)} The fixed field $k(V)^H=k$.

{\rm (d)} All $H$-eigenspaces in $k(V)$ are $1$-dimensional.
\qed\endproclaim

For comparison, we state a ring-theoretic version of this theorem. Since it
will not be needed below, however, we leave the proof to the reader.

\proclaim{4.2\. Theorem} Let $R$ be a commutative affine domain
over $k$, and let $H$ be an algebraic torus acting rationally on $R$ by
$k$-algebra automorphisms. Then the following conditions are equivalent:

{\rm (a)} There are only finitely many $H$-orbits in $\max
R$.

{\rm (b)} Some $H$-orbit is dense in $\max R$.

{\rm (c)} $(\fract R)^H =k$, where $\fract R$ is the quotient field of $R$.

{\rm (d)} All $H$-eigenspaces in $R$ are $1$-dimensional.

{\rm (e)} There are only finitely many $H$-invariant prime
ideals in $R$.
\qed\endproclaim

In the noncommutative setting, consider a prime,
noetherian, affine $k$-algebra $R$, let $\fract R$ be the Goldie
quotient ring of $R$, and replace the maximal ideal space by $\prim
R$. The conditions of Theorem 4.2 are no longer equivalent, since $\prim R$
may well have a dense point (if $0$ is a primitive ideal). It is thus
natural to focus on the stronger conditions, (c) and (d).

 Following Ingalls
\cite{\Ing}, we define a {\it quantum (affine) toric variety\/} to be an
affine domain over
$k$ equipped with a rational action of an algebraic torus $H$ by $k$-algebra
automorphisms, such that the
$H$-eigenspaces are $1$-dimensional. If we fix $H$ and require the action
to be faithful, such algebras can be classified by elements of
$\wedge^2H$ together with finitely generated submonoids that generate
the character group $\widehat H$ \cite{\Ing, Theorem 2.6}.  

Quantum toric
varieties may also be presented as `cocycle twists' (see below) of coordinate
rings of classical toric varieties, analogous to the way quantum tori and
quantum affine spaces are twisted versions of (Laurent) polynomial rings.
Furthermore, quantum toric varieties can be written as factor algebras of
quantum affine spaces, which shows, in particular, that they are noetherian.
It also allows us to apply Theorem 3.2 and express their primitive spectra
as topological quotients of classical varieties. For this purpose, a
presentation as a cocycle twist of a commutative affine algebra suffices;
one-dimensionality of eigenspaces is not needed. Hence, we can work with a
somewhat larger class of algebras than quantum toric varieties.

\definition{4.3} Suppose that $R$ is a $k$-algebra graded by an abelian
group $G$, and that $c: G\times G\rightarrow \kx$ is a $2$-cocycle. The
underlying $G$-graded vector space of $R$ can be equipped with a new
multiplication
$*$, where
$r*s= c(\alpha,\beta)rs$ for all homogeneous elements $r\in R_\alpha$ and
$s\in R_\beta$. This new multiplication is associative, and gives $R$ a new
structure as $G$-graded $k$-algebra, called the {\it
twist of $R$ by $c$\/} (see \cite{\ArScTa,
Section 3} for details). The following alternative presentation is
helpful in keeping computations straight. Let
$R'$ be an isomorphic graded vector space copy of $R$, equipped with a
$G$-graded
$k$-linear isomorphism
$R\rightarrow R'$ denoted $r\mapsto r'$. Then $R'$ becomes a $G$-graded
$k$-algebra such that $r's'=   c(\alpha,\beta)(rs)'$ for $r\in R_\alpha$ and
$s\in R_\beta$. The {\it twist map\/} $r\mapsto r'$ then gives an
isomorphism of the twist of $R$ by $c$ onto $R'$.

Now suppose that $R$ is a commutative affine $k$-algebra, that $G$ is
torsionfree, and that either $-1$ is not in the subgroup of $\kx$ generated
by the image of $c$ or $\chr k=2$. Then there exist topological
quotient maps
$$\spec R\twoheadrightarrow \spec R' \qquad\text{and}\qquad \max
R\twoheadrightarrow \prim R'$$
\cite{\qaffquo, Theorem 6.3}. We shall show that these maps can be given by
formulas analogous to those in Theorem 3.5.

It is convenient to immediately make a reduction to a special choice of
$c$. First, we can replace $G$ by any finitely generated subgroup
containing the support of $R$. Then, as in the proof of \cite{\qaffquo,
Theorem 6.3}, there exists an alternating bicharacter $d$ on $G$,
satisfying the same hypotheses as $c$, such that $R'$ is isomorphic to
the twist of $R$ by $d$. Thus, there is no loss of generality in assuming
that $c$ is an alternating bicharacter, and we shall describe the
topological quotient maps above under this assumption.

The torsionfreeness hypothesis on $G$ is only needed for the reduction
just indicated. Hence, we drop this assumption in the discussion to
follow.
\enddefinition

\definition{4.4} To recap, we are now assuming that $R$ is a commutative
affine $k$-algebra graded by an abelian group $G$, that $c: G\times
G\rightarrow \kx$ is an alternating bicharacter, and that either
$-1\notin \langle \im c\rangle$ or $\chr k=2$. Let $A=R'$ be the twist of
$R$ by $c$.

Choose homogeneous $k$-algebra generators $r_1,\dots,r_n$ for $R$, and
set $a_i= r_i'\in A$ for $i=1,\dots,n$. We use these elements to define sets
$\spec_w R$, $\spec_w A$, $\max_w R$, $\prim_w A$ for $w\subseteq
\{1,\dots,n\}$ along the same lines as in (3.1). Thus
$$\align \spec_w R &= \{P\in\spec R\mid r_i\in P \Longleftrightarrow i\in
w\}\\
\spec_w A &= \{P\in\spec A\mid a_i\in P \Longleftrightarrow i\in
w\}, \endalign$$
while $\max_w R= (\max R)\cap (\spec_w R)$ and $\prim_w A= (\prim A)\cap
(\spec_w A)$. Of course, some of these sets may be empty.

Set $H= \Hom(G,\kx)$, an abelian group under pointwise multiplication.
Because $R$ and $A$ are $G$-graded $k$-algebras, $H$ acts on them by
$k$-algebra automorphisms such that $h{.}r= h(x)r$ and
$h{.}a= h(x)a$ for $h\in H$, $r\in R_x$, $a\in A_x$. Note that the sets
$\spec_w R$ and
$\spec_w A$ are invariant under the induced $H$-actions, since the
elements $r_i$ and $a_i$ are $H$-eigenvectors.

For $i=1,\dots,n$, let $\delta_i\in G$ denote the degree of $r_i$ and
$a_i$. For $w\subseteq
\{1,\dots,n\}$, set $G_w= \sum_{i\notin w} \ZZ\delta_i$ and let $c_w$
denote the restriction of $c$ to $G_w$. By definition of $H$, we have a
pairing $H\times G\rightarrow \kx$, and we use $^\perp$ to denote
orthogonals with respect to this pairing. In particular, for each $w$ we
can define $\rad(c_w)^\perp$, which is a subgroup of $H$.
\enddefinition

We can now state the following more precise version of \cite{\qaffquo,
Theorem 6.3}:

\proclaim{4.5\. Theorem} Let $R$ be a commutative
affine $k$-algebra, graded by an abelian group $G$. Let $c: G\times
G\rightarrow \kx$ be an alternating bicharacter, and let $A=R'$ be the
twist of
$R$ by $c$. Assume that either
$-1\notin \langle \im c\rangle$ or $\chr k=2$. 

With notation as in {\rm (4.4)}, there exist $H$-equivariant topological
quotient maps
$$\phi_s : \spec R\twoheadrightarrow \spec A \qquad\text{and}\qquad \phi_m :
\max R\twoheadrightarrow \prim A$$
such that $P\mapsto (P:\rad(c_w)^\perp)'$ for $P\in \spec_w R$.
Alternatively, $\phi_s(P)$ is the largest ideal of $A$ contained in $P'$.
Moreover, the fibres of $\phi_m$ over points in
$\prim_w A$ are precisely the
$\rad(c_w)^\perp$-orbits in $\max_w R$. \endproclaim

Note that the second description of $\phi_s$ shows that this map depends
only on $R$, $G$, and $c$, that is, it is independent of the choice of
homogeneous generators $r_i\in R$ as in (4.4).

The proof of Theorem 4.5 will be distributed over the following three
subsections.

\definition{4.6} Let $\Phi: A\rightarrow R$ be the inverse of the twist map;
this is a $G$-graded $k$-linear isomorphism. In particular, $X'=
\Phi^{-1}(X)$ for subsets $X\subseteq R$.

Set $\gam= \ZZ^n$ and $\gamp= (\ZZ^+)^n$, and let $\rho: \gam\rightarrow
G$ be the group homomorphism given by the rule $\rho(\alpha_1,
\dots,\alpha_n)= \alpha_1\delta_1 +\dots+ \alpha_n\delta_n$. Set $\ctil=
c\circ (\rho\times\rho)$, which is an alternating bicharacter on $\gam$,
as is $\sigma= \ctil^{\,2}$. For $\alpha= (\alpha_1,
\dots,\alpha_n)\in \gamp$, define
$$r^\alpha= r_1^{\alpha_1}r_2^{\alpha_2} \cdots r_n^{\alpha_n} \in
R_{\rho(\alpha)}
\qquad\text{and}\qquad a_\alpha= (r^\alpha)' \in A_{\rho(\alpha)}.$$
Then $a_\alpha a_\beta= \ctil(\alpha,\beta) a_{\alpha+\beta}$ for
$\alpha, \beta\in \gamp$.

Next, set $\Rtil= \kgamp$ and $\Atil= k^{\ctil}\gamp$, expressed with
$k$-bases $\{y_\alpha\}$ and $\{x_\alpha\}$ such that $y_\alpha y_\beta=
y_{\alpha+\beta}$ and $x_\alpha x_\beta= \ctil(\alpha,\beta)
x_{\alpha+\beta}$ for
$\alpha, \beta\in \gamp$. Let $\Phitil: \Atil\rightarrow \Rtil$ be the
inverse of the twist map; this is a $\gam$-graded $k$-linear isomorphism
such that $\Phitil(x_\alpha)= y_\alpha$ for $\alpha\in \gamp$. There
exist surjective $k$-algebra maps $\tau:
\Rtil\rightarrow R$ and $\pi: \Atil\rightarrow A$ such that
$\tau(y_\alpha)= r^\alpha$ and $\pi(x_\alpha)= a_\alpha$ for $\alpha\in
\gamp$. Thus, we obtain a commutative diagram

\ignore
$$\xymatrixrowsep{3pc}\xymatrixcolsep{5pc}
\xymatrix{
{\Rtil} \ar@{->>}[d]_{\tau} &{\Atil} \ar[l]_{\tilde{\Phi}}
\ar@{->>}[d]^{\pi}\\ 
R &A \ar[l]_{\Phi}
}$$
\endignore

Let $\Htil= \Hom(\gam,\kx)$, acting via $k$-algebra automorphisms on
$\Atil$ and $\Rtil$ in the usual way. For $w\subseteq \{1,\dots,n\}$,
define $\spec_w \Atil$, $\spec_w \Rtil$, $\gam_w$, $\Atil_w$, $\Rtil_w$,
$\sigma_w$, $\ctil_w$ as in (3.1), (3.3), (3.4). (To define $\spec_w \Atil$
and
$\spec_w \Rtil$, for instance, use the elements $x_{\epsilon_i} \in\Atil$
and $y_{\epsilon_i} \in\Rtil$, where $\epsilon_1,\dots,\epsilon_n$ denotes
the standard basis for $\gam$.)  Because of our hypotheses on
$c$, the group
$\langle
\im
\ctil\rangle$ contains no elements of order 2, and hence $\rad(\sigma_w)=
\rad(\ctil_w)$ for all
$w$. We set $\swtil= \rad(\ctil_w)$.

Now Theorem 3.5 gives us $\Htil$-equivariant topological quotient maps
$$\phitil_s: \spec\Rtil \twoheadrightarrow \spec\Atil
\qquad\text{and}\qquad \phitil_m: \max\Rtil \twoheadrightarrow
\prim\Atil$$
such that $P\mapsto \Phitil^{-1}(P:\swtilperp)$ for $P\in
\spec_w \Rtil$. Moreover, the fibres of $\phitil_m$ over points in
$\prim_w \Atil$ are the $\swtilperp$-orbits in $\max_w \Rtil$.

Set $V= \{P\in \spec\Rtil \mid P\supseteq \ker\tau\}$ and $W= \{P\in
\spec\Atil \mid P\supseteq \ker\pi\}$. The proof of \cite{\qaffquo,
Theorem 6.3} shows that $\phitil_s$ and $\phitil_m$ restrict to
topological quotient maps $V\twoheadrightarrow W$ and
$V\cap \max\Rtil \twoheadrightarrow W\cap \prim\Atil$. Therefore there
are topological quotient maps $\phi_s$ and $\phi_m$ such that the
following diagrams commute:

\ignore
$$\xymatrixrowsep{3pc}\xymatrixcolsep{4pc}
\xymatrix{
{\spec\Rtil} \ar@{->>}[r]^{\tilde{\phi}_s} &{\spec\Atil}
&{\max\Rtil} \ar@{->>}[r]^{\tilde{\phi}_m} &{\prim\Atil}\\
V \ar[u]^{\subseteq} \ar@{->>}[r] &W \ar[u]_{\subseteq}
&{V\cap \max\Rtil} \ar[u]^{\subseteq} \ar@{->>}[r] &{W\cap
\prim\Atil} \ar[u]_{\subseteq}\\
{\spec R} \ar[u]^{\tau^{-1}} \ar@{->>}[r]^{\phi_s} &{\spec A}
\ar[u]_{\pi^{-1}} &{\max R} \ar[u]^{\tau^{-1}}
\ar@{->>}[r]^{\phi_m} &{\prim A} \ar[u]_{\pi^{-1}}
}$$
\endignore

\noindent It remains to show that $\phi_s$ and $\phi_m$ have the
properties announced in the statement of Theorem 4.5. 

Given any $P\in \spec R$, we have $\pi^{-1}(\phi_s(P))=
\phitil_s(\tau^{-1}(P))$, which by Lemma 3.6 equals the largest ideal of
$\Atil$ contained in $\Phitil^{-1}\tau^{-1}(P)$. If $I$ is an ideal of $A$
contained in $P'= \Phi^{-1}(P)$, then $\pi^{-1}(I)$ is an ideal of $\Atil$
contained in $\pi^{-1}\Phi^{-1}(P)= \Phitil^{-1}\tau^{-1}(P)$. Hence,
$\pi^{-1}(I) \subseteq \pi^{-1}(\phi_s(P))$, and so $I\subseteq \phi_s(P)$.
This shows that $\phi_s(P)$ is the largest ideal of $A$ contained in $P'$.
\enddefinition

The following lemma
records some facts and conditions that we shall need to finish the proof.

\proclaim{4.7\. Lemma} Let $w\subseteq \{1,\dots,n\}$, and  set $S_w=
\rad(c_w)$.
\roster

\item"(a)" $\ker\tau \subseteq
(\tau^{-1}(P):\swtilperp) \subseteq h(\tau^{-1}(P))$ for $P\in \spec_w R$ and
$h\in \swtilperp$. 
\smallskip

\item"(b)" $(h{.}(-))\circ\tau= \tau\circ ((h\rho){.}(-))$ for $h\in H$.
\smallskip

\item"(c)" $\rho^{-1}(S_w)\cap \gam_w= \swtil$ and so $\rho(\swtil)= S_w$.
\smallskip

\item"(d)"  $(h\rho)(\tau^{-1}(P))= \tau^{-1}(h(P))$ for $P\in \spec_w R$ and
$h\in \swperp$.
\smallskip

\item"(e)" $\swtilperp= \rho^*(\swperp) \gamperp_w$, where $\rho^*:
H\rightarrow \Htil$ denotes the homomorphism given by composition with
$\rho$.
\smallskip

\item"(f)" For $P\in \spec_wR$, the set $\{ \tau^{-1}(h(P)) \mid h\in
\swperp\}$ equals the $\swtilperp$-orbit of $\tau^{-1}(P)$.
\endroster \endproclaim

\demo{Proof} (a) Note that
$\tau^{-1}(P)\in \spec_w \Rtil$ and $\phitil_s(\tau^{-1}(P)) \in W$. Thus
$\Phitil^{-1}(\tau^{-1}(P):\swtilperp) \supseteq \ker\pi$, and so
$(\tau^{-1}(P):\swtilperp) \supseteq \ker\tau$. The remaining inclusion is
clear.

(b) This follows from the observation that
$$h{.}\tau(y_\alpha)= h{.}r^\alpha= h\rho(\alpha)r^\alpha=
\tau(h\rho(\alpha)y_\alpha)= \tau((h\rho){.}y_\alpha)$$
for $h\in H$ and $\alpha\in \gamp$.

(c) Observe that $\rho(\gam_w)= G_w$. Since $S_w\subseteq G_w$ and $\swtil
\subseteq \gam_w$, the second part of the claim will follow from the first.
Now consider an arbitrary element $\gamma\in \gam_w$. Then $\rho(\gamma)\in
S_w$ if and only if $c(\rho(\gamma),-) \equiv 1$ on $G_w= \rho(\gam_w)$, if
and only if $\ctil(\gamma,-) \equiv 1$ on $\gam_w$, if and only if
$\gamma\in \swtil$.

(d) It is clear from part (b) that $\tau(h\rho)(\tau^{-1}(P))= h(P)$.
Since $h\in \swperp$ and $\rho(\swtil)= S_w$ (part (c)), we have $h\rho\in
\swtilperp$, whence $\ker\tau \subseteq (h\rho)(\tau^{-1}(P))$ (part (a)).
Part (d) follows.

(e) Since $\swtil \subseteq \gam_w$, it is clear that $\gamperp_w \subseteq
\swtilperp$. As in the proof of part (d), it is also clear that
$\rho^*(\swperp) \subseteq \swtilperp$.

Now consider $h'\in \swtilperp$; thus $\swtil \subseteq \ker h'$. Hence, $h'$
induces a homomorphism $h_1 : \gam/\swtil \rightarrow \kx$. In view of part
(c), the restriction of $\rho$ to $\gam_w$ induces an isomorphism $\rho_w:
\gam_w/\swtil \rightarrow G_w/S_w$. Let $h_2: G_w\rightarrow \kx$ denote the
composition
$$G_w @>{\text{quo}}>> G_w/S_w @>{\rho_w^{-1}}>> \gam_w/\swtil
@>{\subseteq}>> \gam/\swtil @>{h_1}>> \kx,$$
and observe that the restrictions of $h_2\rho$ and $h'$ to $\gam_w$
coincide. Since $k$ is algebraically closed, $\kx$ is a divisible group, and
hence $\kx$ is injective in the category of abelian groups. Consequently,
$h_2$ can be extended to a homomorphism $h_3\in H$. On one hand, $h_3\rho$
and $h'$ agree on $\gam_w$, whence $h'(h_3\rho)^{-1} \in \gamperp_w$. On the
other hand, $S_w\subseteq \ker h_2\subseteq \ker h_3$, whence $h_3\in
\swperp$. Therefore $h'\in  \rho^*(\swperp) \gamperp_w$, and part (e) is
proved.

(f) If $h\in \swperp$, then $h\rho\in \swtilperp$ and we see from part (d)
that
$\tau^{-1}(h(P))$ lies in the $\swtilperp$-orbit of $\tau^{-1}(P)$.

Conversely, consider $h'\in \swtilperp$. By part (e), $h'= (h\rho)h_0$ for
some $h\in \swperp$ and $h_0\in \gamperp_w$. We claim that $h_0$ leaves
$\tau^{-1}(P)$ invariant; it will then follow from part (d) that
$h'(\tau^{-1}(P))= \tau^{-1}(h(P))$.

Since $P\in \spec_w R$, we see that $\tau^{-1}(P)$ contains the ideal $J_w=
\langle y_{\epsilon_i}
\mid i\in w\rangle$, where $\epsilon_1,\dots,\epsilon_n$
denotes the standard basis for
$\gam$. Note that $J_w$ is invariant under $H$. Since
$\Rtil/J_w$ is spanned by the cosets $y_\alpha+J_w$ for $\alpha\in
\gam_w\cap \gamp$, we see that the induced action of $\gamperp_w$ on
$\Rtil/J_w$ is trivial. Thus $\tau^{-1}(P)$ is indeed invariant under $h_0$,
which establishes part (f). \qed\enddemo

\definition{4.8} (Proof of Theorem 4.5) Let $w\subseteq \{1,\dots,n\}$,
and  consider $P\in \spec_w R$. In view of Lemma 4.7(a) and the
commutative diagrams in (4.6), we see that
$$\align \phi_s(P) &= \pi\phitil_s(\tau^{-1}(P))=
\pi\Phitil^{-1}(\tau^{-1}(P):\swtilperp)\\
 &= \Phi^{-1}\tau(\tau^{-1}(P):\swtilperp)= \bigl(
\tau(\tau^{-1}(P):\swtilperp)
\bigr)'. \endalign$$
On the other hand, Lemma 4.7(f) implies that
$$(\tau^{-1}(P):\swtilperp)= \bigcap_{h\in\swperp} \tau^{-1}(h(P))=
\tau^{-1}(P:\swperp),$$
whence $\tau(\tau^{-1}(P):\swtilperp)= (P:\swperp)$. Therefore $\phi_s(P)=
(P:\swperp)'$.

Since $H$ is abelian, the set functions $(-:\swperp)$ are $H$-equivariant,
as is the twist map $(-)'$. Hence, it follows from the
formula just proved that $\phi_s$ is $H$-equivariant.

Finally, note that the fibre of $\phi_m$ over any point of $\prim_w A$ is
contained in
$\max_w R$. Consider $P,P'\in \max_w R$. Taking account of (4.6), we see
that $\phi_m(P)=
\phi_m(P')$ if and only if $\phitil_m(\tau^{-1}(P))=
\phitil_m(\tau^{-1}(P'))$, if and only if $\tau^{-1}(P)$ and
$\tau^{-1}(P')$ lie in the same $\swtilperp$-orbit. By Lemma 4.7(f), this
occurs if and only if $\tau^{-1}(P')= \tau^{-1}(h(P))$ for some $h\in
\swperp$, and thus if and only if $P$ and $P'$ lie in the same
$\swperp$-orbit. \qed\enddefinition

Once the precise form of the maps in Theorem 4.5 is
given, we can easily construct topological quotient maps between different
cocycle twists of $R$, as follows.

\proclaim{4.9\. Corollary} Let $R$ be a commutative
affine $k$-algebra, graded by an abelian group $G$. Let $c_1,c_2 : G\times
G\rightarrow \kx$ be alternating bicharacters, and let $A_i$ be the
twist of
$R$ by $c_i$. Assume that either
$-1\notin \langle \im c_1\rangle \cup \langle \im c_2\rangle$ or $\chr k=2$.
Set $H=\Hom(G,\kx)$, and let $\phi_i : \spec R\twoheadrightarrow \spec A_i$
be the
$H$-equivariant topological quotient map given in Theorem 4.5.

Assume that $c_2(x,y)=1 \implies c_1(x,y)=1$, for any $x,y\in G$. Then there
exists an $H$-equi\-var\-i\-ant topological quotient map $\phi: \spec A_1
\twoheadrightarrow \spec A_2$ such that the following diagram commutes:

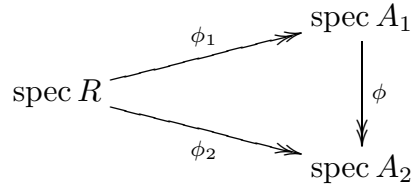
\ignore
$$\xymatrixrowsep{0.75pc}\xymatrixcolsep{5pc}
\xymatrix{
 &\spec A_1 \ar@{->>}[dd]^{\phi}\\
\spec R \ar@{->>}[ur]^{\phi_1} \ar@{->>}[dr]_{\phi_2}\\
 &\spec A_2
}$$
\endignore \endproclaim

\demo{Proof} Set up notation for $A_i$ and $\phi_i$ as in (4.4), and
abbreviate $(c_i)_w$ by $c_{iw}$. Let $\Phi_i: A_i\rightarrow R$ be the
inverse of the twist map. Then $\phi_i$ is given by the rule
$$\phi_i(P)= \Phi_i^{-1}(P: \rad(c_{iw})^\perp) \qquad \text{for\ \ } P\in
\spec_w R.$$
Our hypotheses on $c_1,c_2$ imply that $\rad(c_{2w}) \subseteq \rad(c_{1w})$
for all $w$, whence $\rad(c_{2w})^\perp \supseteq \rad(c_{1w})^\perp$ and so
$$(P: \rad(c_{2w})^\perp)= \bigl( (P: \rad(c_{1w})^\perp):
\rad(c_{2w})^\perp \bigr)$$
for all $P\in \spec_w R$. Hence, we can define a map $\phi : \spec A_1
\rightarrow \spec A_2$ such that
$$\phi(Q)= \Phi_2^{-1}( \Phi_1(Q): \rad(c_{2w})^\perp)$$
for $Q\in \spec_w A_1$. It is clear that $\phi\phi_1 =\phi_2$. Since
$\phi_1$ and $\phi_2$ are $H$-equivariant topological quotient maps, so is
$\phi$. \qed\enddemo

\enddocument